\documentclass[9pt,shortpaper,twoside,web]{IEEEtran}
\usepackage{}
\ifCLASSINFOpdf
\else
\fi

\usepackage{graphicx}
\usepackage{dcolumn}
\usepackage{bm}

\usepackage{graphics} 
\usepackage{epsfig} 
\usepackage{times} 
\usepackage{amsmath} 
\usepackage{amssymb}  
\usepackage{accents}
\usepackage{amsfonts}
\usepackage{fancyhdr}
\usepackage{color}
\usepackage{subfigure}

\newtheorem{theorem}{Theorem}
\newtheorem{lemma}[theorem]{Lemma}

\newtheorem{proposition}[theorem]{Proposition}

\newtheorem{definition}{Definition}

\newtheorem{remark}{Remark}

\allowdisplaybreaks[4] 

\begin{document}
%
\title{A {De~Giorgi} Iteration-based Approach for the Establishment of ISS Properties for Burgers' Equation with Boundary and In-domain Disturbances}
%
%
%

\author{Jun~Zheng$^{1}$\thanks{$^{1}$School of Civil Engineering and School of Mathematics, Southwest Jiaotong University, Chengdu, Sichuan, P. R. of China 611756
       {\tt\small zhengjun2014@aliyun.com}}and Guchuan~Zhu$^{2}$,~\IEEEmembership{Senior~Member,~IEEE}\thanks{$^{2}$Department of Electrical Engineering, Polytechnique Montr\'{e}al,  P.O. Box 6079, Station Centre-Ville, Montreal, QC, Canada H3T 1J4
       {\tt\small guchuan.zhu@polymtl.ca}}

\thanks{\textcolor{blue}{This paper has been accepted for publication by IEEE TAC, and is available at http://dx.doi.org/10.1109/TAC.2018.2880160}}

}

%
%

\markboth{Manuscript Submitted to IEEE Trans. on Automatic Control}%
{Zheng \MakeLowercase{\textit{et al.}}: }
%




\maketitle

\begin{abstract}
This note addresses input-to-state stability (ISS) properties with respect to (w.r.t.) boundary and in-domain disturbances for Burgers' equation. The developed approach is a combination of the method of De~Giorgi iteration and the technique of Lyapunov functionals by adequately splitting the original problem into two subsystems. The ISS properties in $L^2$-norm for Burgers' equation have been established using this method. Moreover, as an application of De~Giorgi iteration, ISS in $L^\infty$-norm w.r.t. in-domain disturbances and actuation errors in boundary feedback control for a 1-$D$ {linear} {unstable reaction-diffusion equation} have also been established. It is the first time that the method of De~Giorgi iteration is introduced in the ISS theory for infinite dimensional systems, and the developed method can be generalized for tackling some problems on multidimensional spatial domains and to a wider class of nonlinear {partial differential equations (PDEs)}.
\end{abstract}


%
\IEEEpeerreviewmaketitle

\section{Introduction}\label{Sec: Introduction}
Extending the theory of ISS, which was originally developed for finite-dimensional nonlinear systems \cite{Sontag:1989,Sontag:1990}, to infinite dimensional systems has received a considerable attention in the recent literature. In particular, there are significant progresses on the establishment of ISS estimates with respect to disturbances \cite{Argomedo:2013, Argomedo:2012, Dashkovskiy:2010, Dashkovskiy:2013, jacob2016input, Karafyllis:2014, Karafyllis:2016, Karafyllis:2016a, Karafyllis:2018, Logemann:2013, Mazenc:2011, Mironchenko:2018, Prieur:2012, Tanwani:2017,Zheng:2017} for different types of {PDEs}.

It is noticed that most of the earlier work on this topic dealt with disturbances distributed over the domain. It was demonstrated that the method of Lyapunov functionals is a well-suited tool for dealing with a wide rang of problems of this category. Moreover, it is shown in \cite{Argomedo:2012} that the method of Lyapunov functionals can be readily applied to some systems with boundary disturbances by transforming the later ones to a distributed disturbance. {However, ISS estimates obtained by such a method may include time derivatives of boundary disturbances, which is not strictly in the original form of ISS formulation.

The problems with disturbances acting on the boundaries usually lead to a formulation involving unbounded operators. It is shown in \cite{jacob2016input,Jacob:2017} that for a class of linear PDEs, the exponential stability plus a certain admissibility implies the ISS and iISS (integral input-to-state stability \cite{jacob2016input,Sontag:1998}) w.r.t. boundary disturbances. However, it may be difficult to assess this property for nonlinear PDEs. To resolve this concern while not invoking unbounded operators in the analysis, it is proposed in \cite{Karafyllis:2016,Karafyllis:2016a,karafyllis2017siam} to derive the ISS property directly from the estimates of the solution to the considered PDEs using the method of spectral decomposition and finite-difference. ISS in $L^2$-norm and in weighted $L^\infty$-norm for PDEs with a Sturm-Liouville operator is established by applied this method in \cite{Karafyllis:2016,Karafyllis:2016a,karafyllis2017siam}. However, spectral decomposition and finite-difference schemes may involve heavy computations for nonlinear PDEs or problems on multidimensional spatial domains. It is introduced in \cite{Mironchenko:2017} a monotonicity-based method for studying the ISS of nonlinear parabolic equations with boundary disturbances. It is shown that with the monotonicity the ISS of the original nonlinear parabolic PDE with constant boundary disturbances is equivalent to the ISS of a closely related nonlinear parabolic PDE with constant distributed disturbances and zero boundary conditions. As an application of this method, the ISS properties in $L^p$-norm ({$\forall p>2$}) for some linear parabolic systems have been established.
In a recent work \cite{Zheng:2017}, the classical method of Lyapunov functionals is applied to establish ISS properties in $L^2$-norm w.r.t. boundary disturbances for a class of semilinear parabolic PDEs. Some technical inequalities have been developed, which allows dealing directly with items on the boundary points in proceeding on ISS estimates. {}{The result of \cite{Zheng:2017}} shows that the method of Lyapunov functionals is still effective in obtaining ISS properties for some linear and nonlinear PDEs with Neumann or Robin boundary conditions. However, the technique used may not be suitable for problems with Dirichlet boundary conditions.

{}{The present work is dedicated to the establishment of ISS properties for Burgers' equation that is one of the most popular PDEs in mathematical physics \cite{Hopf:1991}. Burgers' equation is considered as a simplified form of the Navier-Stokes equation and can be used to approximate the Saint-Venant equation. Therefore, the study on the control of the Burgers' equation is an important and natural step for flow control and many other fluid dynamics inspired applications. Indeed, there exists a big amount of work on the control of the Burgers' equation, e.g., just to cite a few, \cite{Azmi:2016_SIAM,Burns:1991, Byrnes:1998, Krstic:1999, Krstic:2008_IEEE_TAC, Liu:2000, Liu:2001}.}

{}{The problem dealt with in this work can be seen as a complementary setting compared to that considered in \cite{Zheng:2017} in the sense that the problem is subject to Dirichlet boundary conditions.} The method developed in this note consists first in splitting the original system into two subsystems: one system with boundary disturbances and a zero-initial condition, and another one with no boundary disturbances, but with homogenous boundary conditions and a non-zero initial condition. Note that the in-domain disturbances can be placed in either of these two subsystems. Then, ISS properties in $L^\infty$-norm for the first system will be deduced by the technique of De~Giorgi iteration, and ISS properties in $L^2$-norm (or $L^\infty$-norm) for the second system will be established by the method of Lyapunov functionals. Finally, the ISS properties in $L^2$-norm (or $L^\infty$-norm) for the original system are obtained by combining the ISS properties of the two subsystems. With this method, we established the ISS in $L^2$-norm for Burgers' equation with boundary and in-domain disturbances. Moreover, using the techniques of transformation, splitting and De~Giorgi iteration, we established the ISS in $L^\infty$-norm for a 1-$D$ {}{linear unstable reaction-diffusion equation} with boundary feedback control including actuation errors. Note that although the De~Giorgi iteration is a classic method in regularity analysis of elliptic and parabolic PDEs, it is the first time, to the best of our knowledge, that it is introduced in the investigation of ISS properties for PDEs. Moreover, the technique of De~Giorgi iteration may be applicable for certain {}{nonlinear PDEs and }problems on multidimensional spatial domains.

The rest of the note is {}{organized} as follows. Section~\ref{Sec: Preliminaries} introduces briefly the technique of De~Giorgi iteration and presents some preliminary inequalities needed for the subsequent development. Section~\ref{Sec: Main results} presents the considered problems and the main results. Detailed development on the establishment of ISS properties for Burgers' equation is given Section~\ref{Sec: Burgers Eq}. The application of De~Giorgi iteration in the establishment of ISS in  $L^\infty$-norm for a 1-$D$ {}{linear unstable reaction-diffusion equation} is illustrated in Section~\ref{Sec: reaction-diffusion Eq}. Finally, some concluding remarks are provided in Section~\ref{Sec: Conclusion}.

\section{Preliminaries}\label{Sec: Preliminaries}

\subsection{{De~Giorgi} iteration}\label{Sec: De Giorgi iteration}

De~Giorgi iteration is an important tool for regularity analysis of elliptic and parabolic PDEs. In his famous work on linear elliptic equations published in 1957 \cite{DeGiorgi:1957}, De~Giorgi established local boundedness and H\"{o}lder continuity for functions satisfying certain integral inequalities, known as the De~Giorgi class of functions, which completed the solution of Hilbert's 19$^{\text{th}}$ problem. The same problem has been resolved independently by Nash in 1958 \cite{Nash1958}. It was shown later by Moser that the result of De~Giorgi and Nash can be obtained using a different formulation \cite{Moser1961}. In the literature, this method is often called the De~Giorgi-Nash-Moser theory.

{}{Let $\mathbb{R}:=(-\infty,+\infty)$}, $\Omega\subset \mathbb{R}^n(n\geq 1)$ be an open bounded set, and $\gamma $ be a constant. The De~Giorgi class $DG^+(\Omega,\gamma)$ consists of functions $u\in W^{1,2}(\Omega)$ which satisfy, for every ball $B_r(y)\subset \Omega$, every $0<r'<r$, and every $k \in \mathbb{R}$, the following Caccioppoli type inequality:
\begin{align*}
\int_{B_{r'}(y)}|\nabla (u-k)_+|^2\text{d}x\leq \frac{\gamma}{(r-r')^2}\int_{B_{r}(y)}| (u-k)_+|^2\text{d}x,
\end{align*}
where $ (u-k)_+=\max\{u-k,0\}$. The class $DG^-(\Omega,\gamma)$ is defined in a similar way. The main idea of De~Giorgi iteration is to estimate $ |A_k|$, the measure of $\{x\in \Omega;u(x)\geq k\}$, and derive $|A_k|=0$ with some $k$ for functions $u$ in De~Giorgi class by using the iteration formula given below.
\begin{lemma}[{}{\cite[Lemma 4.1.1]{Wu2006}}]\label{iteration}
Suppose that $\varphi$ is a non-negative decreasing function on $[k_0,+\infty)$ satisfying
\begin{align*}
\varphi(h)\leq \bigg(\frac{M}{h-k}\bigg)^\alpha\varphi^\beta(k),\ \ \forall h>k\geq k_0,
\end{align*}
where $M>0,\alpha>0,\beta>1$ are constants. Then the following holds
\begin{align*}
\varphi(k_0+l_0)=0,
\end{align*}
with $l_0=2^{\frac{\beta}{\beta-1}}M(\varphi(k_0))^{\frac{\beta-1}{\alpha}}$.
\end{lemma}

The method of De~Giorgi iteration can be generalized to some linear parabolic PDEs and PDEs with a divergence form (see, e.g., \cite{DiBenedetto:2010,Wu2006}). However, this method in its original formulation cannot be applied directly in the establishment of ISS properties for infinite dimensional systems. The main reason is that the obtained boundedness of a solution depends always on {some data that is increasing in $t$ rather than a class $\mathcal {K}\mathcal {L}$ function associated with $u_0 $ and $t$}, even for linear parabolic PDEs \cite{DiBenedetto:2010,Wu2006}, {which is not under the form of ISS}. To overcome this difficulty, we developed in this work an approach that amounts first to splitting the original problem into two subsystems and then to applying the De~Giorgi iteration together with the technique of Laypunov functionals to obtain the ISS estimates of the solutions expressed in the standard formulation of the ISS theory.

\subsection{Preliminary inequalities}\label{Sec: preliminary results}
{}{Let $\mathbb{R}_+:=(0,+\infty)$ and $\mathbb{R}_{\geq 0} :=  [0,+\infty)$. For notational simplicity, we always denote $\|\cdot\|_{L^{2}(0,1)}$ by $\|\cdot\|$ in this note.} We present below two inequalities needed for the subsequent development.
\newline
\begin{lemma}\label{Lemma 3}
Suppose that $u\in C^{1}([a,b];\mathbb{R})$, then for any $p\geq 1$, {}{one has}
\begin{align}\label{eq: Sobolev embedding}
\bigg(\int_{a}^b|u|^p\text{d}x\bigg)^{\frac{1}{p}}\leq (b-a)^{\frac{1}{p}}\bigg(\frac{2}{b-a}\|u\|^2+(b-a)\|u_x\|^2\bigg)^{\frac{1}{2}}.
\end{align}
\end{lemma}
\begin{IEEEproof}
We show first that
\begin{align}
u^2(c)\leq \dfrac{2}{b-a}\|u\|^2+(b-a)\|u_x\|^2,\ \forall c\in[a,b].\label{001}
\end{align}
Denote $(u_z(z))^2$ by $u_z^2(z)$. For each $c\in [a,b]$, let $g(x)={\int_{c}^xu^2_z(z)\text{d}z}$.  Note that $g'(x)=u^2_x(x)$. By H\"{o}lder's inequality (see \cite[Appendix B.2.e]{Evans:2010}), we have
\begin{align*}
\displaystyle\left(\int_{x}^cu_z(z)\text{d}z\right)^2 \leq \left|(c-x)\int_{x}^cu^2_z(z)\text{d}z\right|
 =(x-c)g(x).
\end{align*}
It follows
\begin{align*}
u^2(c)&=\bigg(u(x)+{\int_{x}^cu_z(z)\text{d}z}\bigg)^2\notag\\
&\leq 2u^2(x)+2\bigg({\int_{x}^cu_z(z)\text{d}z}\bigg)^2\notag\\
&\leq 2u^2(x)+2(x-c)g(x).
\end{align*}
Integrating over $[a,b]$ and
noting that
\begin{align*}
&\int_{a}^b(x-c) g(x)\text{d}x\notag\\
=&\bigg[\frac{(x-c)^2}{2}g(x)\bigg]\bigg|_{x=a}^{x=b}-\int_{a}^b\frac{(x-c)^2}{2} u^2_x(x) \text{d}x\notag\\
\leq &\frac{(b-c)^2}{2}\int_{c}^bu^2_x(x)\text{d}x- \frac{(a-c)^2}{2}\int_{c}^au^2_x(x)\text{d}x\notag\\
\leq &\frac{(b-a)^2}{2}\int_{a}^bu^2_x(x)\text{d}x,
\end{align*}
we get
$
u^2(c)(b-a)\leq 2 \|u\|^2 + (b-a)^2\|u_x\|^2,
$
which yields \eqref{001}.

Now by \eqref{001}, we have
\begin{align*}
\bigg(\int_{a}^b|u|^p\text{d}{x}\bigg)^{\frac{1}{p}}&\leq \bigg(\int_{a}^b\max_{x\in[a,b]}|u|^p\text{d}{x}\bigg)^{\frac{1}{p}}\\
&=  (b-a)^{\frac{1}{p}}\max_{x\in[a,b]}|u|\\
&\leq
(b-a)^{\frac{1}{p}}\bigg(\frac{2}{b-a}\|u\|^2+(b-a)\|u_x\|^2\bigg)^{\frac{1}{2}}.
\end{align*}
\end{IEEEproof}
\begin{remark}
{}{Note first that \eqref{eq: Sobolev embedding} is a variation of Sobolev embedding inequality, which will be used in the De~Giorgi iteration in the analysis of the Burgers' equation with in-domain and Dirichlet boundary disturbances. Moreover, the inequality~\eqref{001} is an essential technical result for the establishment of the ISS w.r.t. boundary disturbances for PDEs with Robin or Neumann boundary conditions (see, e.g., \cite{Zheng:2017}). Therefore, these two inequalities play an important role in the establishment of the ISS for PDEs w.r.t. boundary and in-domain disturbances with either Robin, or Neumann, or Dirichlet boundary conditions.}
\end{remark}

\section{Problem Formulation and Main Results}\label{Sec: Main results}

\subsection{Problem formulation and well-posedness analysis}\label{Sec: Problem formulations}
In this work, we address ISS properties for Burgers' equation with Dirichlet boundary conditions:
\begin{subequations}\label{++28}
\begin{align}
&u_t-\mu u_{xx}+\nu uu_x=f(x,t)\ \ {\text{in}\  (0,1)\times\mathbb{R}_{+}},\label{++28'}\\
 &u(0,t)=0,u(1,t)=d(t),\label{++2}\\
  &u(x,0)=u_0(x),\label{++3}
\end{align}
\end{subequations}
where $\mu>0$, $\nu>0$ are constants, $d(t)$ is the disturbance on the boundary, which can represent actuation and sensing errors, and the function $f(x,t)$ is the disturbance distributed over the domain. {}{Throughout this note, we always assume that $f\in \mathcal {H}^{\theta,\frac{\theta}{2}}([0,1]\times \mathbb{R}_{\geq 0})$ and $d\in \mathcal {H}^{1+\frac{\theta}{2}}(\mathbb{R}_{\geq 0})$ for some $\theta\in (0,1)$.}

{}{We refer to \cite[Chapter~1, pages 7-9]{Ladyzenskaja:1968} for the definition on H\"{o}lder type function spaces $ \mathcal {H}^{l}([0,1])$, $ \mathcal {H}^{l}([0,T])$, $\mathcal {H}^{l,\frac{l}{2}}([0,1]\times[0,T])$, $C^1([0,1])$, $C^{2,1}([0,1]\times[0,T])$, $\mathcal {H}^{l}(\mathbb{R}_{\geq 0})$, $\mathcal {H}^{l,\frac{l}{2}}([0,1]\times \mathbb{R}_{\geq 0})$ and $C^{2,1}([0,1]\times\mathbb{R}_{\geq 0})$, where $l>0$ is a nonintegral number and $T>0$.} {}{We also refer to \cite[Chapter~1, page~12]{Ladyzenskaja:1968} for the statement of classical solutions of Cauchy problems.}


{}{The result for well-posedness assessment of \eqref{++28} is given below, which is guaranteed by \cite[Theorem 6.1, pages 452-453]{Ladyzenskaja:1968}.
\begin{proposition}\label{Proposition 1}
Assume that $u_0\in \mathcal {H}^{2+\theta}([0,1])$ with {$u_0(0)=0$, $u_0(1)=d(0)$, $\mu u_0{''}(0)+f(0,0)=0$ and $\mu u_0{''}(1)+f(1,0) =d{'}(0)$}. For any $T>0$, there exists a unique classical solution $u\in \mathcal {H}^{2+\theta,1+\frac{\theta}{2}}( [0,1]\times [0,T])\subset C^{2,1}( [0,1]\times[0,T])$ of \eqref{++28}.
\end{proposition}}
\begin{remark}
The proof of Proposition~\ref{Proposition 1} follows from Theorem 6.1 in \cite[pages 452-453]{Ladyzenskaja:1968}, which establishes the existence of a unique solution in the H\"{o}lder space of functions {$ \mathcal {H}^{2+\theta,1+\frac{\theta}{2}}( [0,1]\times [0,T])$} for a more general quasilinear parabolic equations with Dirichlet boundary conditions. It should be noticed that the proof of Theorem~6.1 in \cite[pages 452-453]{Ladyzenskaja:1968} is based on the linearization of the considered system and the application of the Leray-Schauder theorem on fixed points. Since {$\mathcal {H}^{2+\theta,1+\frac{\theta}{2}}( [0,1]\times [0,T])\subset C^{2,1}( [0,1]\times[0,T])$}, we can obtain the existence of the unique classical solution in the time interval $[0,T]$, where $T>0$ can be arbitrarily large.
\end{remark}

\subsection{{Main results on ISS estimates for Burgers' equation}}
Let $\mathcal {K}=\{\gamma : \mathbb{R}_{\geq 0} \rightarrow \mathbb{R}_{\geq 0}|\ \gamma(0)=0,\gamma$ is continuous, strictly increasing$\}$; $ \mathcal {K}_{\infty}=\{\theta \in \mathcal {K}|\ \lim\limits_{s\rightarrow\infty}\theta(s)=\infty\}$; $ \mathcal {L}=\{\gamma : \mathbb{R}_{\geq 0}\rightarrow \mathbb{R}_{\geq 0}|\ \gamma$ is continuous, strictly decreasing, $\lim\limits_{s\rightarrow\infty}\gamma(s)=0\}$; $ \mathcal {K}\mathcal {L}=\{\beta : \mathbb{R}_{\geq 0}\times \mathbb{R}_{\geq 0}\rightarrow \mathbb{R}_{\geq 0}|\ \beta(\cdot,t)\in \mathcal {K}, \forall t \in \mathbb{R}_{\geq 0}$, and $\beta(s,\cdot)\in \mathcal {L}, \forall s \in {}{\mathbb{R}_{+}}\}$.

{\begin{definition}
System~\eqref{++28} is said to be input-to-state stable (ISS) in $L^q$-norm ($2\leq q\leq +\infty$) {w.r.t. {boundary disturbances} $d\in \mathcal {H}^{1+\frac{\theta}{2}}(\mathbb{R}_{\geq 0})$ and {in-domain disturbances} $f\in \mathcal {H}^{\theta,\frac{\theta}{2}}([0,1]\times \mathbb{R}_{\geq 0})$}, if there exist functions $\beta\in \mathcal {K}\mathcal {L}$ and $ \gamma_1, \gamma_2,\in \mathcal {K}$ such that the solution to \eqref{++28} satisfies
\begin{align}\label{Eq: ISS def2}
\begin{split}
\|{}{u(\cdot,t)}\|_{L^{q}(0,1)}\leq & \beta\left( \|{u_0}\|_{L^{q}(0,1)},t\right)+\gamma_1\left(\max_{s\in [0,t]}|d(s)|\right)\\
      &+\gamma_2\left(\max_{(x,s)\in [0,1]\times [0,t]}|f(x,s)|\right),\ \forall t\geq 0.
\end{split}
\end{align}
System~\eqref{++28} is said to be ISS {w.r.t. boundary disturbances $d\in \mathcal {H}^{1+\frac{\theta}{2}}(\mathbb{R}_{\geq 0})$, and integral input-to-state stable (iISS) w.r.t. in-domain disturbances $f\in \mathcal {H}^{\theta,\frac{\theta}{2}}([0,1]\times \mathbb{R}_{\geq 0})$}, in $L^q$-norm ($2\leq q\leq +\infty$),
if there exist functions $\beta\in \mathcal {K}\mathcal {L},\theta\in \mathcal {K}_{\infty} $, and $\gamma_1 ,\gamma_2 \in \mathcal {K}$ such that the solution to \eqref{++28} satisfies
\begin{align}\label{Eq: iISS def2}
\begin{split}
\|{}{u(\cdot,t)}\|_{L^{q}(0,1)} \leq & \beta\left( \|{u_0}\|_{L^{q}(0,1)},t\right)
      +\gamma_1\left(\max_{s\in [0,t]}|d(s)|\right) \\
      &+\theta\left(\!\!\int_{0}^t\!\!\gamma_2(\|f(\cdot,s)\|)\text{d}s\right),\ \forall t\geq 0.
\end{split}
\end{align}
Moreover, System~\eqref{++28} is said to be exponential input-to-state stable (EISS), or exponential integral input-to-state stable (EiISS),  w.r.t. {boundary disturbances $d(t)$, or in-domain disturbances $f(x,t)$}, if there exist $\beta'\in \mathcal {K}_{\infty}$ and a constat $\lambda > 0$ such that  $\beta( \|{u_0}\|_{L^{q}(0,1)},t) \leq \beta'(\|{u_0}\|_{L^{q}(0,1)})e^{-\lambda t}$ in \eqref{Eq: ISS def2} or \eqref{Eq: iISS def2}.
\end{definition}}

In order to apply the technique of splitting and the method of De~Giorgi iteration in the investigation of the ISS properties for the considered problem, {while guaranteeing the well-posedness by Proposition~\ref{Proposition 1} for every system}, we assume that the compatibility condition {$u_0(0)=u_0''(0)=u_0(1)=u_0''(1)=d(0)=d'(0)=f(0,0)=f(1,0)=0$} always holds {}{in Section \ref{Sec: Main results} and \ref{Sec: Burgers Eq}}. Furthermore, unless} stated, we always take a certain function in {}{$C^{2,1}( [0,1]\times\mathbb{R}_{\geq 0})$} as the unique solution of a considered system. Then the ISS properties w.r.t. boundary and in-domain disturbances for System~\eqref{++28} are stated in the following {theorems}.

\begin{theorem} \label{Theorem 11}
{System \eqref{++28} is {}{EISS} in $L^2$-norm w.r.t. {boundary disturbances} $d\in \mathcal {H}^{1+\frac{\theta}{2}}(\mathbb{R}_{\geq 0})$ and {in-domain disturbances} $f\in \mathcal {H}^{\theta,\frac{\theta}{2}}([0,1]\times \mathbb{R}_{\geq 0})$ satisfying $\sup\limits_{s\in \mathbb{R}_{\geq 0}} |d(s)| + \frac{4\sqrt{2}}{\mu}\sup\limits_{(x,s)\in[0,1]\times \mathbb{R}_{\geq 0}}|f(x,s)| < \frac{\mu}{\nu}$}, having the following estimate for any $t>0$:
 \begin{align*}
\|u(\cdot,t)\|^2\leq & 2\|u_0\|^2 e^{-\mu t} +4\max\limits_{s\in [0,t)} |d(s)|^2 \nonumber \\
                     &+\frac{128}{\mu^2}\max\limits_{(x,s)\in[0,1]\times [0,t]}|f(x,s)|^2.
\end{align*}
\end{theorem}
\begin{theorem} \label{Theorem 11-2}
{System \eqref{++28} is {}{EISS} in $L^2$-norm w.r.t. {boundary disturbances} $d\in \mathcal {H}^{1+\frac{\theta}{2}}(\mathbb{R}_{\geq 0})$ satisfying $\sup\limits_{t\in \mathbb{R}_{\geq 0}} |d(t)|<\frac{\mu}{\nu}$, and {}{EiISS} w.r.t. {in-domain disturbances} $f\in \mathcal {H}^{\theta,\frac{\theta}{2}}([0,1]\times \mathbb{R}_{\geq 0})$,} having the following estimate for any $t>0$:
 \begin{align*}
\|u(\cdot,t)\|^2 \leq& 2\|u_0\|^2 e^{-(\mu -\varepsilon)t}+2\max\limits_{s\in [0,t]}|d(s)|^2 \nonumber \\
                     &+\frac{2}{\varepsilon} \int_{0}^t\|f(\cdot,s)\|^2\text{d}s,\ \forall\varepsilon\in (0,\mu).
\end{align*}
\end{theorem}
\begin{remark}
In general, the boundness of the disturbances is a reasonable assumption for nonlinear PDEs in the establishment of ISS properties \cite{Mironchenko:2016}. However, as shown in Section \ref{Sec: reaction-diffusion Eq}, the boundedness of the disturbances may be not a necessary condition for ISS properties of linear PDEs.
\end{remark}
{}{\begin{remark}
{}{As pointed out in \cite{Karafyllis:2016a}, the} assumptions on the continuity of $f$ and $d$ are required for assessing the existence of a classical solution of the considered system. However, they are only sufficient conditions and can be weakened {}{if solutions in a weak sense are considered. Moreover, for the establishment of ISS estimates, the assumptions on the continuity of $f$ and $d$ can eventually be relaxed.}
\end{remark}}

\section{Proofs of ISS Estimates for Burgers' Equation}\label{Sec: Burgers Eq}
\subsection{Proof of Theorem~\ref{Theorem 11}}
In this section, we establish the ISS estimates for Burgers' equation w.r.t. boundary and in-domain disturbances described in Theorem~\ref{Theorem 11} by using the technique of splitting. Specifically, let $w$ be the unique solution of the following system:
\begin{subequations}\label{++29}
 \begin{align}
 &w_t-\mu w_{xx}+\nu ww_x=f(x,t)\ \ \text{in}\ (0,1)\times\mathbb{R}_{+},\\
 &w(0,t)=0,w(1,t)=d(t),\\
  &w(x,0)=0.
\end{align}
\end{subequations}
Then $v=u-w$ is the unique solution of the following system:
\begin{subequations}\label{++31}
 \begin{align}
 &v_t-\mu v_{xx}+\nu vv_x+\nu {(wv )_x}=0 \ \text{in}\ (0,1)\times\mathbb{R}_{+}, \\
 &v(0,t)=v(1,t)=0,\\
 &v(x,0)=u_0(x).
\end{align}
\end{subequations}

For System~\eqref{++29}, we have the following estimate.
\begin{lemma} \label{Theorem 12} Suppose that $\mu>0,\nu>0$. For every $t>0$, {}{one has}
\begin{align}\label{++13}
&\max\limits_{(x,s)\in[0,1]\times [0,t]} |w(x,s)| \notag\\
&\;\;\; \leq \max\limits_{s\in [0,t]}|d(s)| + \frac{4\sqrt{2}}{\mu} \max\limits_{(x,s)\in[0,1]\times [0,t]}|f(x,s)|.
\end{align}
\end{lemma}
For System~\eqref{++31}, we have the following estimate.
\begin{lemma} \label{Theorem 13}{Suppose that $\mu>0,\nu>0$, and ${\sup\limits_{t\in \mathbb{R}_{\geq 0}}} |d(t)|+ \frac{4\sqrt{2}}{\mu} {\sup\limits_{(x,t)\in [0,1]\times\mathbb{R}_{\geq 0}}} |f(x,t)|<\frac{\mu}{\nu}$}. For every $t>0$, {}{one has}
\begin{align*}
\|v(\cdot,t)\|^2\leq {\|u_0\|^2 e^{-\mu t}.}
\end{align*}
\end{lemma}
Then the result of Theorem~\ref{Theorem  11} is a consequence of {Lemma}~\ref{Theorem  12} and {Lemma}~\ref{Theorem  13}.
\begin{IEEEproof}[Proof of Theorem~\ref{Theorem  11}] Note that $u=w+v$, we get by {Lemma}~\ref{Theorem  12} and {Lemma}~\ref{Theorem  13}:
\begin{align*}
\|u(\cdot,t)\|^2
\leq & 2\|w(\cdot,t)\|^2+2\|v(\cdot,t)\|^2\\
\leq & 2\left(\max\limits_{(x,s)\in [0,1]\times[0,t]} |w(x,s)|\right)^2+2\|v(\cdot,t)\|^2\\
\leq & 2\|u_0\|^2 e^{-\mu t}\notag\\
&+2\left(\!\!\max\limits_{s\in [0,t]}|d(s)|+\frac{4\sqrt{2}}{\mu}\!\! \max\limits_{(x,s)\in [0,1]\times[0,t]} \!\!|f(x,s)|\!\!\right)^2.
\end{align*}
\end{IEEEproof}
In the following, we use De~Giorgi iteration and Lyapunov method to prove {Lemma}~\ref{Theorem 12} and {Lemma}~\ref{Theorem 13}, respectively.

\begin{IEEEproof}[Proof of Lemma~\ref{Theorem 12}] In order to apply the technique of De~Giorgi iteration, we shall define some quantities.
For any $t>0$, let $k_0=\max\Big\{\max\limits_{s\in[0,t]}d(s),0\Big\}$. For any $k\geq k_0$, let $ \eta(x,s)=(w(x,s)-k)_+\chi_{[t_1,t_2]}(
s)$, where $\chi_{[t_1,t_2]}(s) $ is the character function on $[t_1,t_2]$ and $0\leq t_1<t_2\leq t$. Let $A_{k}(s)=\{x\in (0,1);w(x,s)>k\}$ and {}{$\varphi_{k}=\sup\limits_{s\in(0,t)}|A_{k}(s)|$}, where $|B|$ denotes the 1-dimensional Lebesgue measure of a set $B\subset(0,1)$. For any $p>2$, let $l_0= \frac{1}{\mu}2^{\frac{5p-8}{2p-4}}\max\limits_{(x,s)\in [0,1]\times[0,t]}|f(x,s)| \varphi_{k_0}$. The main idea of De~Giorgi iteration is to show that $|A_{k_0+l_0}(s)|=0 $ for almost every $s\in [0,t]$, which yields $~~~~{\text{ess}}\!\!\!\!\!\!\!\!\sup\limits_{\!\!\!\!\!\!\!\!(x,s)\in[0,1]\times [0,t]} w(x,s)\leq k_0+l_0 $. The lower boundedness of $w(x,s)$ can be obtained in a similar way. Then the desired result is guaranteed by the continuity of $w$ and its lower and upper boundedness.

{}{Although the computation of De~Giorgi iteration can follow a standard process (see, e.g., the case of linear parabolic equations presented in \cite[Theorem 4.2.1, \S4.2.2]{Wu2006}), we provide the details for completeness.} Multiplying \eqref{++29} by $\eta$, and {}{ noting that $(w(0,s)-k)_+ =(w(1,s)-k)_+=0$ for $k\geq k_0$ and $s\in [0,t]$}, we get
\begin{align}\label{+16}
&\int_{0}^t\int_{0}^1(w-k)_t(w-k)_+\chi_{[t_1,t_2]}(s)\text{d}x\text{d}s \notag\\
&+\mu\int_{0}^t\int_{0}^1|((w-k)_+)_x |^2\chi_{[t_1,t_2]}(s)\text{d}x\text{d}s\notag\\
&+\nu\int_{0}^t\int_{0}^1ww_x(w-k)_+\chi_{[t_1,t_2]}(s) \text{d}x\text{d}s \notag\\
=&\int_{0}^t\int_{0}^1f(w-k)_+\chi_{[t_1,t_2]}(s) \text{d}x\text{d}s.
\end{align}
Let $I_k(s)=\int_{0}^1((w(x,s)-k)_+)^2\text{d}x$, which is absolutely continuous on $[0,t]$. Suppose that $I_k(t_0)=\max\limits_{s\in[0,t]}I_k(s)$ with some $t_0\in [0,t]$. Due to $I_k(0)=0$ and $I_k(s)\geq 0 $, we can assume that $t_0>0$ without loss of generality.

For $ \varepsilon>0$ small enough, choosing $t_1=t_0-\varepsilon$ and $t_2=t_0 $, it follows
\begin{align*}
&\frac{1}{2\varepsilon}\int_{t_0-\varepsilon}^{t_0}\frac{d}{dt}\int_{0}^1((w-k)_+)^2\text{d}x\text{d}s \notag\\
&+\frac{\mu}{\varepsilon}\int_{t_0-\varepsilon}^{t_0}\int_{0}^1|((w-k)_+)_x |^2\text{d}x\text{d}s \notag\\
&+\frac{\nu}{\varepsilon}\int_{t_0-\varepsilon}^{t_0}\int_{0}^1ww_x(w-k)_+ \text{d}x\text{d}s \notag\\
\leq & \frac{1}{\varepsilon}\int_{t_0-\varepsilon}^{t_0}\int_{0}^1|f|(w-k)_+\text{d}x\text{d}s.
\end{align*}
Note that
 \begin{align*}
\frac{1}{2\varepsilon}\int_{t_0-\varepsilon}^{t_0}\frac{d}{dt}
\int_{0}^1((w-k)_+)^2\text{d}x\text{d}s&=\frac{1}{2\varepsilon}( I_k(t_0)-I_k(t_0-\varepsilon))\notag\\
&\geq 0.
\end{align*}
We have
\begin{align*}
     & \frac{\mu}{\varepsilon}\int_{t_0-\varepsilon}^{t_0}\int_{0}^1|((w-k)_+)_x |^2\text{d}x\text{d}s \notag\\
     &+\frac{\nu}{\varepsilon}\int_{t_0-\varepsilon}^{t_0}\int_{0}^1ww_x(w-k)_+ \text{d}x\text{d}s \notag\\
\leq & \frac{1}{\varepsilon}\int_{t_0-\varepsilon}^{t_0}\int_{0}^1|f|(w-k)_+\text{d}x\text{d}s.
\end{align*}
Letting $ \varepsilon\rightarrow 0^+$, and noting that
\begin{align}
&\lim_{\varepsilon\rightarrow 0^+}\frac{1}{\varepsilon}\int_{t_0-\varepsilon}^{t_0}\int_{0}^1ww_x(w-k)_+ \text{d}x\text{d}s \notag\\
= &\int_{0}^1w(x,t_0)w_x(x,t_0)(w(x,t_0)-k)_+ \text{d}x \notag\\
=& \int_{0}^1(w(x,t_0)-k)_+((w(x,t_0)-k)_+)_x(w(x,t_0)-k)_+ \text{d}x \notag\\
 &+\int_{0}^1k((w(x,t_0)-k)_+)_x(w(x,t_0)-k)_+ \text{d}x\notag\\
=& \frac{1}{3}((w(x,t_0)-k)_+)^{3}|^{x=1}_{x=0}+\frac{k}{2}((w(x,t_0)-k)_+)^{2}|^{x=1}_{x=0}\notag\\
=&0,\label{+201803}
\end{align}
we get
\begin{align}
&\mu\int_{0}^1|((w(x,t_0)-k)_+)_x |^2\text{d}x \notag\\
\leq & \int_{0}^1|f(x,t_0)|(w(x,t_0)-k)_+\text{d}x.\label{+20181}
\end{align}
We deduce by Lemma~\ref{Lemma 3}, Poincar\'{e}'s inequality \cite[Chap.~2, Remark~2.2]{Krstic:2008}, and \eqref{+20181} that for any $p>2$,
\begin{align*}
     &\bigg(\int_{0}^1|(w(x,t_0)-k)_+ |^p\text{d}x\bigg)^{\frac{2}{p}} \notag\\
\leq & 2 \int_{0}^1|((w(x,t_0)-k)_+)_x |^2\text{d}x \notag\\
\leq & \frac{2}{\mu} \int_{0}^1|f(x,t_0)|(w(x,t_0)-k)_+\text{d}x.
\end{align*}

 Then we have
\begin{align*}
&\bigg(\int_{A_{k}(t_0)}|(w(x,t_0)-k)_+ |^p\text{d}x\bigg)^{\frac{2}{p}} \notag\\
\leq & \frac{2}{\mu}   \int_{A_{k}(t_0)}|f(x,t_0)|(w(x,t_0)-k)_+\text{d}x.
\end{align*}
By H\"{o}lder's inequality (see \cite[Appendix B.2.e]{Evans:2010}), it follows
\begin{align*}
&\bigg(\int_{A_{k}(t_0)}|(w(x,t_0)-k)_+ |^p\text{d}x\bigg)^{\frac{2}{p}} \notag\\
\leq &\frac{2}{\mu} \bigg(\int_{A_{k}(t_0)}|(w(x,t_0)-k)_+|^p\text{d}x\bigg)^{\frac{1}{p}}\bigg(\int_{0}^1|f(x,t_0)|^q\text{d}x\bigg)^{\frac{1}{q}},
\end{align*}
where $\frac{1}{p}+\frac{1}{q}=1$.
Thus
\begin{align}\label{++14}
&\bigg(\int_{A_{k}(t_0)}|(w(x,t_0)-k)_+ |^p\text{d}x\bigg)^{\frac{1}{p}}\notag\\
\leq & \frac{2}{\mu} \bigg(\int_{A_{k}(t_0)}|f(x,t_0)|^q\text{d}x\bigg)^{\frac{1}{q}}\notag\\
\leq & \frac{2}{\mu} |{A_{k}(t_0)}|^{\frac{1}{q}} \max_{(x,s)\in [0,1]\times[0,t]}|f(x,s)|\notag\\
\leq & \frac{2}{\mu} \max_{(x,s)\in [0,1]\times[0,t]}|f(x,s)|\varphi_{k}^{\frac{1}{q}}.
\end{align}
Now for $I_k(t_0)$, we get by H\"{o}lder's inequality and \eqref{++14}
\begin{align*}
I_k(t_0)&\leq \bigg(\int_{A_{k}(t_0)}|(w(x,t_0)-k)_+ |^p\text{d}x \bigg)^{\frac{2}{p}}|{A_{k}(t_0)}|^{\frac{p-2}{p}}\notag\\
&\leq \bigg(\frac{2}{\mu} \max_{(x,s)\in [0,1]\times[0,t]}|f(x,s)|\bigg)^2\varphi_{k}^{3-\frac{4}{p}}.
\end{align*}
Recalling the definition of $I_k(t_0)$, for any $s\in [0,t]$ we conclude that
\begin{align}
{I_k(s)}\leq I_k(t_0)\leq \bigg(\frac{2}{\mu} \max_{(x,s)\in [0,1]\times[0,t]}|f(x,s)|\bigg)^2\varphi_{k}^{3-\frac{4}{p}}.\label{++15}
\end{align}
Note that for any $h>k$ and $s\in [0,t]$ the following holds
\begin{align}
{I_k(s)}\geq \int_{A_{h}({}{s})}\!\!\!\!\!\!\!\!((w(x,{}{s})-k)_+)^2 \text{d}x\geq (h-k)^2{|A_h(s)|}.\label{+201802}
\end{align}
Then we infer from \eqref{++15} and \eqref{+201802} that
\begin{align*}
(h-k)^2\varphi_h\leq \bigg(\frac{2}{\mu} \max_{(x,s)\in [0,1]\times[0,t]}|f(x,s)|\bigg)^2\varphi_{k}^{3-\frac{4}{p}},
\end{align*}
which is
\begin{align*}
\varphi_h\leq \left(\frac{2}{\mu}\frac{\max\limits_{(x,s)\in [0,1]\times[0,t]}|f(x,s)|}{h-k}\right)^2\varphi_{k}^{3-\frac{4}{p}}.
\end{align*}
As $p>2$, we have $ 3-\frac{4}{p}>1$. By Lemma~\ref{iteration}, we obtain
\begin{align*}
\varphi_{k_0+l_0}=\sup_{s\in[0,t]}|A_{k_0+l_0}|=0,
\end{align*}
where $l_0=2^{\frac{3p-4}{2p-4}}\frac{2}{\mu} \max\limits_{(x,s)\in [0,1]\times[0,t]}|f(x,s)|\varphi_{k_0}^{1-\frac{2}{p}}\leq \frac{1}{\mu}2^{\frac{5p-8}{2p-4}}\max\limits_{(x,s)\in [0,1]\times[0,t]}|f(x,s)|$.

By the definition of $A_k$, for almost {every} $(x,s)\in [0,1]\times [0,t]$, one has
\begin{align*}
w(x,s)
\leq & k_0+\frac{1}{\mu}2^{\frac{5p-8}{2p-4}}\max\limits_{(x,s)\in [0,1]\times[0,t]}|f(x,s)|\notag\\
 = &\max\Big\{\max\limits_{s\in[0,t]}d(s),0\Big\} +\frac{1}{\mu}2^{\frac{5p-8}{2p-4}}\max\limits_{(x,s)\in [0,1]\times[0,t]}|f(x,s)|.
\end{align*}
By continuity of $w(x,s) $, for every $(x,s)\in [0,1]\times [0,t]$, the following holds
\begin{align*}
w(x,s)
\leq \max\Big\{\max\limits_{s\in[0,t]}d(s),0\Big\}\frac{1}{\mu}2^{\frac{5p-8}{2p-4}}\max\limits_{(x,s)\in [0,1]\times[0,t]}|f(x,s)|.
\end{align*}
Letting $p\rightarrow +\infty$, we get for every $(x,s)\in [0,1]\times [0,t]$
\begin{align}
w(x,s)
\leq & \max\Big\{\max\limits_{s\in[0,t]}d(s),0\Big\}\notag\\
     &+\frac{4\sqrt{2}}{\mu} \max\limits_{(x,s)\in [0,1]\times[0,t]}|f(x,s)|.\label{++16}
\end{align}
To conclude on the inequality \eqref{++13}, we need also to prove the lower boundedness of $w(x,t)$. Indeed, setting $\overline{w}=-w$, we get
\begin{align*}
 &\overline{w}_t-\mu \overline{w}_{xx}-\nu \overline{w}\overline{w}_x=-f(x,t),\\
 &\overline{w}(0,t)=0,\overline{w}(1,t)=-d(t),\\
 &\overline{w}(x,0)=0.\notag
\end{align*}
 Proceeding as above and noting \eqref{+201803}, the following equality holds in the process of De~Giorgi iteration:
\begin{align*}
\lim_{\varepsilon\rightarrow 0^+}\frac{1}{\varepsilon}\int_{t_0-\varepsilon}^{t_0}\int_{0}^1-\overline{w}\overline{w}_x(\overline{w}-k)_+ \text{d}x\text{d}s =0.
\end{align*}
Then for every $(x,s)\in [0,1]\times [0,t]$ we have
\begin{align}
- w(x,s)=&\overline{w}(x,s) \notag\\
           \leq & \max\Big\{\max\limits_{s\in[0,t]}-d(s),0\Big\}
                 +\frac{4\sqrt{2}}{\mu} \max\limits_{(x,s)\in [0,1]\times[0,t]}|f(x,s)|.\label{++17''}
\end{align}
Finally, \eqref{++13} follows from \eqref{++16} and \eqref{++17''}.
\end{IEEEproof}
\begin{IEEEproof}[Proof of {Lemma}~\ref{Theorem 13}]
Multiplying \eqref{++31} by $v$ and integrating over $(0,1)$, we get
\begin{align*}
\int_{0}^1\!\!v_tv\text{d}x+\mu\int_{0}^1v^2_{x}\text{d}x+\nu\int_{0}^1v^2v_x\text{d}x+\nu\int_{0}^1(wv)_xv\text{d}x =0.
\end{align*}
Note that $\int_{0}^1v^2v_x\text{d}x=\frac{1}{3}v^3|^{x=1}_{x=0}=0$ and
\begin{align*}
\int_{0}^1(wv)_xv\text{d}x = wv^2 |^{x=1}_{x=0}-\int_{0}^1wvv_x\text{d}x=-\int_{0}^1wvv_x\text{d}x.
\end{align*}
By Young's inequality (see \cite[Appendix B.2.d]{Evans:2010}), H\"{o}lder's inequality (see \cite[Appendix B.2.e]{Evans:2010}), {Lemma}~\ref{Theorem 12}, and {the assumption on $d$}, we deduce that
\begin{align}
&\frac{1}{2}\frac{d}{dt}\|v(\cdot,t)\|^2+\mu\|v_x(\cdot,t)\|^2 \leq \nu\int_{0}^1|wvv_x|\text{d}x\notag\\
\leq &\frac{\nu}{2}\max\limits_{(x,s)\in[0,1]\times [0,t]} |w(x,s)|(\|v(\cdot,t)\|^2+\|v_x(\cdot,t)\|^2) \notag\\
\leq &\frac{\nu}{2}\bigg(\max\limits_{s\in [0,t]}|d(s)|+\frac{4\sqrt{2}}{\mu} \max\limits_{(x,s)\in [0,1]\times[0,t]} |f(x,s)|\bigg)\notag\\
&\times(\|v(\cdot,t)\|^2+\|v_x(\cdot,t)\|^2)\notag\\
\leq &\frac{\nu}{2}\times\frac{\mu}{\nu}(\|v(\cdot,t)\|^2+\|v_x(\cdot,t)\|^2)
\notag\\
=&\frac{\mu}{2}(\|v(\cdot,t)\|^2+\|v_x(\cdot,t)\|^2).\label{+201804}
\end{align}

By Poincar\'{e}'s inequality, we have
\begin{align}
\mu\|v_x(\cdot,t)\|^2&=\frac{\mu}{2}\|v_x(\cdot,t)\|^2
+\frac{\mu}{2}\|v_x(\cdot,t)\|^2\notag\\
&\geq \frac{\mu}{2}\|v_x(\cdot,t)\|^2+\mu\|v(\cdot,t)\|^2.\label{+201805}
\end{align}

Then by \eqref{+201804} and \eqref{+201805}, it follows
\begin{align*}
\frac{\text{d}}{\text{d}t}\|v(\cdot,t)\|^2
\leq-\mu\|v(\cdot,t)\|^2,
\end{align*}
which together with Gronwall's inequality (\cite[Appendix B.2.j]{Evans:2010}) yields
\begin{align*}
\|v(\cdot,t)\|^2
&\leq\|v(\cdot,0)\|^2 e^{-\mu t}=\|u_0\|^2 e^{-\mu t}.
\end{align*}
\end{IEEEproof}


\subsection{Proof of Theorem \ref{Theorem 11-2}}
In order to prove Theorem~\ref{Theorem 11-2}, we consider the following two systems:
\begin{subequations}\label{+++32}
 \begin{align}
 &w_t-\mu w_{xx}+\nu ww_x=0\ \ \text{in}\  (0,1)\times\mathbb{R}_{+},\\
 &w(0,t)=0,w(1,t)=d(t),\\
  &w(x,0)=0,
\end{align}
\end{subequations}
and
\begin{subequations}\label{+++33}
 \begin{align}
 &v_t-\mu v_{xx}+\nu vv_x+\nu (wv + vw)_x=f(x,t)\ \ \text{in}\ (0,1)\times\mathbb{R}_{+},\\
 &v(0,t)=v(1,t)=0,\\
  &v(x,0)=u_0(x),
\end{align}
\end{subequations}
where $v=u-w$.

For System~\eqref{+++32}, we have the following estimate, which is a special case (i.e. $f(x,t)=0$) of {Lemma}~\ref{Theorem 12}.
\begin{lemma} \label{Theorem 14} Suppose that $\mu>0,\nu>0$. For every $t>0$, {}{one has}
\begin{align}\label{++32''}
\max\limits_{(x,s)\in[0,1]\times [0,t]} |w(x,s)|\leq \max\limits_{s\in [0,t]}|d(s)|.
\end{align}
\end{lemma}

For System~\eqref{+++33}, we have the following estimate.
\begin{lemma} \label{Theorem 15}{Suppose that $\mu>0,\nu>0$, and ${\sup\limits_{t\in \mathbb{R}_{\geq 0}}} |d(t)|<\frac{\mu}{\nu}$}. For every $t>0$, {}{one has}
\begin{align*}
\|v(\cdot,t)\|^2\leq \|u_0\|^2 e^{-(\mu -\varepsilon)t}+\frac{1}{\varepsilon} \int_{0}^t\|f(\cdot,s)\|^2\text{d}s,\ \forall \varepsilon\in (0,\mu).
\end{align*}
\end{lemma}
Note that $u=w+v$. Then the result of Theorem \ref{Theorem 11-2} is a consequence of {Lemma}~\ref{Theorem 14} and {Lemma}~\ref{Theorem 15}, which can be proven as in Theorem~\ref{Theorem 11}.
\begin{IEEEproof}[Proof of {Lemma}~\ref{Theorem 15}] Multiplying \eqref{+++33} by $v$ and integrating over $(0,1)$, we get
\begin{align*}
&\int_{0}^1v_tv\text{d}x+\mu\int_{0}^1v^2_{x}\text{d}x+\nu\int_{0}^1v^2v_x\text{d}x+\nu\int_{0}^1(wv)_xv\text{d}x \\
=&\int_{0}^1f(x,t)v\text{d}x.
\end{align*}
{}{Arguing as in \eqref{+201804},} we get
\begin{align}
&\frac{1}{2}\frac{d}{dt}\|v(\cdot,t)\|^2+\mu\|v_x(\cdot,t)\|^2\notag \\
\leq &\nu\int_{0}^1|wvv_x|\text{d}x+\int_{0}^1f(x,t)v\text{d}x\notag \\
\leq &\frac{\nu}{2}\max\limits_{s\in [0,t]}|d(s)|(\|v(\cdot,t)\|^2+\|v_x(\cdot,t)\|^2)\notag\\
&+\frac{1}{2\varepsilon}\|f(\cdot,t)\|^2
+\frac{\varepsilon}{2}\|v(\cdot,t)\|^2\notag \\
\leq & \frac{\nu}{2}\frac{\mu}{\nu}(\|v(\cdot,t)\|^2+\|v_x(\cdot,t)\|^2)
+\frac{1}{2\varepsilon}\|f(\cdot,t)\|^2+\frac{\varepsilon}{2}\|v(\cdot,t)\|^2\notag \\
=&\frac{1}{2}(\varepsilon+\mu)\|v(\cdot,t)\|^2
+\frac{\mu}{2}\|v_x(\cdot,t)\|^2+\frac{1}{2\varepsilon}\|f(\cdot,t)\|^2,\label{+201806}
\end{align}
where we choose $0<\varepsilon <\mu$.

By \eqref{+201805} and \eqref{+201806}, we get
\begin{align*}
\frac{d}{dt}\|v(\cdot,t)\|^2
\leq-(\mu -\varepsilon)\|v(\cdot,t)\|^2+\frac{1}{\varepsilon}\|f(\cdot,t)\|^2.
\end{align*}
By Growall's inequality (see \cite[Appendix B.2.j]{Evans:2010}), we have
\begin{align*}
\|v(\cdot,t)\|^2
&\leq\|v(\cdot,0)\|^2 e^{-(\mu -\varepsilon)t}+\frac{1}{\varepsilon} \int_{0}^t\|f(\cdot,s)\|^2\text{d}s\\
&=\|u_0\|^2 e^{-(\mu -\varepsilon)t}+\frac{1}{\varepsilon} \int_{0}^t\|f(\cdot,s)\|^2\text{d}s.
\end{align*}
\end{IEEEproof}

\section{Application to a 1-$D$ {}{Linear} Unstable Reaction-Diffusion Equation with Boundary Feedback Control}\label{Sec: reaction-diffusion Eq}
In this section, we {}{illustrate the application of the developed method in the study of the ISS property for the following 1-$D$ {}{linear} reaction-diffusion equation with an unstable term}:
\begin{align}\label{++9181}
 u_t-\mu u_{xx}+a(x)u=f(x,t) \ \ \ \ \text{in}\ \  (0,1)\times\mathbb{R}_{+},
\end{align}
where $\mu >0$ is a constant, $a\in C^1([0,1])$ and $f\in \mathcal {H}^{\theta,\frac{\theta}{2}}([0,1]\times \mathbb{R}_{\geq 0})$. The system is subject to the boundary and initial conditions
\begin{subequations}\label{++9182}
\begin{align}
 &u(0,t)=0,u(1,t)=U(t),\\
 &u(x,0)=u_0(x),
\end{align}
\end{subequations}
where $U(t)\in \mathbb{R}$ is the control input. {}{Note that the control input can be placed on either ends of the boundary. Nevertheless, it can be switched to the other end by a spatial variable transformation $x\rightarrow 1-x$. The ISS properties of this system w.r.t. boundary disturbances, i.e., $f(x,t)\equiv 0$, have been addressed in \cite{Karafyllis:2016,Karafyllis:2016a,Mironchenko:2017}.}

The stabilization of \eqref{++9181} in a disturbance-free setting with $\mu =1$ and $f(x,t)\equiv 0$ is presented in \cite{Krstic:2008,Liu:2003,Smyshlyaev:2004}. {}{The exponential stability is achieved by means of a backstepping boundary feedback control of the form}
\begin{align}
U(t)=\int_{0}^1k(1,y)u(y,t)\text{d}y,\ \forall t\geq 0,\label{++9183}
\end{align}
where $k\in C^2([0,1]\times[0,1])$ can be obtained as the Volterra kernel of a Volterra integral transformation
\begin{align}
w(x,t)=u(x,t)-\int_{0}^xk(x,y)u(y,t)\text{d}y,\label{+2018032701}
\end{align}
which transforms \eqref{++9181}, \eqref{++9182}, and \eqref{++9183} to the problem \begin{align*}
 w_t-\mu w_{xx}+\nu w=0 \ \ \ \ \text{in}\ \  (0,1)\times\mathbb{R}_{+},
\end{align*}
with $\nu> 0$, subject to the boundary and initial conditions
\begin{align*}
&w(0,t)=w(1,t)=0,\\
&w(x,0)=w_0(x)=u_0(x)-\int_{0}^xk(x,y)u_0(y)\text{d}y.
\end{align*}

When $\mu >0$, $f\in \mathcal {H}^{\theta,\frac{\theta}{2}}([0,1]\times \mathbb{R}_{\geq 0})$, and in the presence of actuation errors represented by the disturbance $d\in \mathcal {H}^{1+\frac{\theta}{2}}(\mathbb{R}_{\geq 0})$, the applied control action is of the form \cite{Karafyllis:2016,Karafyllis:2016a,Mironchenko:2017}
\begin{align}
U(t)=d(t)+\int_{0}^1k(1,y)u(y,t)\text{d}y,\ \forall t\geq 0.\label{++9188}
\end{align}
We can use the Volterra integral transformation \eqref{+2018032701} to transform \eqref{++9181}, \eqref{++9182}, and \eqref{++9188} to the following system
\begin{align}\label{++9185}
 w_t-\mu w_{xx}+\nu w=f(x,t) \ \ \ \ \text{in}\ \  (0,1)\times\mathbb{R}_{+},
\end{align}
with $\nu>0$, subject to the boundary and initial conditions
\begin{subequations}\label{++9189}
\begin{align}
&w(0,t)=0,w(1,t)=d(t),\\
&w(x,0)=w_0(x)=u_0(x)-\int_{0}^xk(x,y)u_0(y)\text{d}y.
\end{align}
\end{subequations}
Then the solution to \eqref{++9181}, \eqref{++9182}, and \eqref{++9188} can be found by the
inverse Volterra integral transformation
\begin{align}
u(x,t)=w(x,t)+\int_{0}^xl(x,y)w(y,t)\text{d}y,\label{++9187}
\end{align}
where $l\in C^2([0,1]\times[0,1])$ is an appropriate kernel. Indeed, the existence of
the kernels $k\in C^2([0,1]\times[0,1])$ and $l\in C^2([0,1]\times[0,1])$ can be obtained in the same way as in \cite{Liu:2003,Smyshlyaev:2004}.

For the system~\eqref{++9181} with \eqref{++9182} and \eqref{++9188}, we have the following ISS estimate.
\begin{proposition}\label{Theorem 16}
Suppose that $\mu >0$, $a\in C^1([0,1])$, {$d\in \mathcal {H}^{1+\frac{\theta}{2}}(\mathbb{R}_{\geq 0})$, $ f\in \mathcal {H}^{\theta,\frac{\theta}{2}}([0,1]\times \mathbb{R}_{\geq 0})$}, and $u_0\in \mathcal {H}^{2+\theta}([0,1])$ for some $\theta\in (0,1)$, with the compatibility conditions:
\begin{align*}
&{u_0(0)=d(0)=d'(0)=f(0,0)=f(1,0)=0,}\\
&u_0(1)=\int_{0}^1k(1,y)u_0(y)\text{d}y,\\
&{u_0(1)\frac{dk(x,x)}{dx}\bigg|_{x=1}+u_0'(1)k(1,1)=0.}
\end{align*}
{System} \eqref{++9181} with \eqref{++9182} and \eqref{++9188} is EISS in $L^\infty$-norm {w.r.t. {boundary disturbances} $d\in \mathcal {H}^{1+\frac{\theta}{2}}(\mathbb{R}_{\geq 0})$ and {in-domain disturbances} $f\in \mathcal {H}^{\theta,\frac{\theta}{2}}([0,1]\times \mathbb{R}_{\geq 0})$}, having the following estimate:
\begin{align*}
 \max_{x\in[0,1]}|u(x,t)|\leq  &C_0 \max_{x\in [0,1]}|u_0|e^{-\nu t}+C_1\bigg(\max\limits_{s\in [0,t]}|d(s)|\notag\\ &+ \frac{4\sqrt{2}}{\mu} \max\limits_{(x,s)\in[0,1]\times [0,t]}|f(x,s)|\bigg),
\end{align*}
where $\nu>0$ is the same as in \eqref{++9185}, $C_1= \Big(1+\max\limits_{(x,y)\in[0,1]\times[0,1] }|l(x,y)|\Big)$ and $C_0= C_1\Big(1+ \max\limits_{(x,y)\in[0,1]\times[0,1] }|k(x,y)|\Big)$ are positive constants.
\end{proposition}
\begin{IEEEproof}
{Note that by the compatibility conditions, it follows $w_0(0)=w_0(1)=w_0''(0)=w_0''(1)=f(0,0)=f(1,0)=d(0)=d'(0)=0 $. Therefore, we can use the technique of splitting and De~Giorgi iteration to establish the ISS estimate for System~\eqref{++9185} with \eqref{++9189}. Let} $g$ be the unique solution of the following system:
\begin{subequations}\label{subequ.1}
\begin{align}
&g_t-\mu g_{xx}+\nu g=f(x,t),\ \ \ \ (0,1)\times\mathbb{R}_{+},\\
&g(0,t)=0,g(1,t)=d(t),\\
&g(x,0)=0,
\end{align}
\end{subequations}
and let $h=w-g$ be the unique solution of the following system:
\begin{subequations}\label{subequ.2}
\begin{align}
&h_t-\mu h_{xx}+\nu h=0,\ \ \ \ (0,1)\times\mathbb{R}_{+},\\
&h(0,t)=h(1,t)=0,\\
&h(x,0)=h_0(x)=w_0(x).
\end{align}
\end{subequations}
For \eqref{subequ.1} and \eqref{subequ.2}, we claim that for any $t\in\mathbb{R}_{\geq 0}$:
\begin{align}\label{estimate1}
&\max\limits_{(x,s)\in[0,1]\times [0,t]} |g(x,s)| \notag\\
&\;\;\; \leq \max\limits_{s\in [0,t]}|d(s)| + \frac{4\sqrt{2}}{\mu} \max\limits_{(x,s)\in[0,1]\times [0,t]}|f(x,s)|,
\end{align}
and
 \begin{align}\label{estimate2}
\max_{x\in [0,1]}|h(x,t)|\leq \max_{x\in [0,1]}|h_0(x)| e^{- \nu t}.
\end{align}
We prove  \eqref{estimate1} by De~Giorgi iteration. Indeed, for any fixed $t>0$, letting $k_0,k$, $ \eta(x,s)$, and $t_0$ be defined as in the proof of Theorem~\ref{Theorem 11} (replace $w$ by $h$) and taking $\eta(x,s)$ as a test function, we get
\begin{align*}
&\int_{0}^t\int_{0}^1(g-k)_t(g-k)_+\chi_{[t_1,t_2]}(s)\text{d}x\text{d}s \notag\\
&+\mu\int_{0}^t\int_{0}^1\chi_{[t_1,t_2]}(s)|((g-k)_+)_x |^2\text{d}x\text{d}s\notag\\
&+\nu\int_{0}^t\int_{0}^1g(g-k)_+\chi_{[t_1,t_2]}(s) \text{d}x\text{d}s \notag\\
=&\int_{0}^t\int_{0}^1f(g-k)_+\chi_{[t_1,t_2]}(s) \text{d}x\text{d}s.
\end{align*}
Noting that $\nu\int_{0}^t\int_{0}^1g(g-k)_+\chi_{[t_1,t_2]}(s) \text{d}x\text{d}s\geq 0$, it can be seen that \eqref{+20181} still holds (replace $w$ by $h$), which leads to \eqref{estimate1}.

For the proof of \eqref{estimate2}, we choose the following Lyapunov functional
\begin{align*}
E(t)=\int_{0}^{1}|h(\cdot,t)|^{2p}\text{d}x,\ \forall p\geq 1,\forall t\in\mathbb{R}_{\geq 0}.
\end{align*}
Applying Poincar\'{e}'s inequality, it follows $\|h^{p}(\cdot,t)\|^2\leq \frac{p^2}{2} \|h^{p-1}h_x(\cdot,t)\|^2$. Then by direct computations, we get
\begin{align*}
\frac{\text{d}}{\text{d}t}E(t)\leq -2p\bigg( \nu +\frac{2\mu (2p-1)}{p^2}\bigg)E(t),
\end{align*}
{}{which together with Gronwall's inequality yields}
\begin{align}\label{estimate2'}
\|h(\cdot,t)\|_{L^{2p}(0,1)}^{2p}\leq \|h_0\|_{L^{2p}(0,1)}^{2p} e^{-2p\big( \nu +\frac{2\mu (2p-1)}{p^2}\big)t},\ \forall p\geq 1.
\end{align}
Taking the $2p$-th root of \eqref{estimate2'} and letting $p\rightarrow +\infty$, it follows
\begin{align}
\|h(\cdot,t)\|_{L^{\infty}(0,1)}\leq \|h_0\|_{L^{\infty}(0,1)} e^{- \nu t},\ \forall t\in\mathbb{R}_{\geq 0} .\label{estimate2''}
\end{align}
Finally, we obtain \eqref{estimate2} by \eqref{estimate2''} and the continuity of $h$ and $h_0$.

As a consequence of \eqref{estimate1} and \eqref{estimate2}, the following estimate holds for any $t\in\mathbb{R}_{\geq 0}$:
\begin{align}\label{estimate3}
\max\limits_{x\in[0,1]} |w(x,t)|\leq &\max\limits_{(x,s)\in[0,1]\times [0,t]} |g(x,s)|+\max_{x\in [0,1]}|h(x,t)| \notag\\
\leq &\max\limits_{s\in [0,t]}|d(s)| + \frac{4\sqrt{2}}{\mu} \max\limits_{(x,s)\in[0,1]\times [0,t]}|f(x,s)|\notag\\
&+\max_{x\in [0,1]}|w_0(x)| e^{- \nu t}.
\end{align}
Note that
\begin{align}
 \max_{x\in [0,1]}|w_0(x)|&\leq \max_{x\in [0,1]}\bigg|u_0-\int_{0}^xk(x,y)u_0(y)\text{d}y \bigg|\notag\\
 &\leq \max_{x\in [0,1]}|u_0|+\max_{x\in [0,1]}\bigg|\int_{0}^xk(x,y)u_0(y)\text{d}y \bigg|\notag\\
 &\leq \bigg(1+ \max\limits_{(x,y)\in[0,1]\times[0,1] }|k(x,y)|\bigg)\max_{x\in [0,1]}|u_0|.\label{estimate4}
\end{align}
{}{Finally, the desired result follows from \eqref{++9187}, \eqref{estimate3}, and \eqref{estimate4}.}
\end{IEEEproof}
\begin{remark}
If we put $f(x,t)$ in \eqref{subequ.2} instead of in \eqref{subequ.1}, and proceed as in the proof of Theorem~\ref{Theorem 11-2}, we can prove that the system \eqref{++9181} with \eqref{++9182} and \eqref{++9188} is EISS w.r.t. {boundary disturbances} and EiISS w.r.t. {in-domain disturbances}.
\end{remark}
\begin{remark}
In the case where $a(x)\equiv a$ is a constant, the ISS in $L^2$-norm and $L^p$-norm ($\forall p>2$) for the system \eqref{++9181} with \eqref{++9182} w.r.t. actuation errors for boundary feedback control \eqref{++9188} is established in \cite{Karafyllis:2016a} by the technique of eigenfunction expansion, and in \cite{Mironchenko:2017} by the monotonicity method, respectively.
\end{remark}
\begin{remark}
The ISS in a weighted $L^\infty$-norm w.r.t. boundary and in-domain disturbances for solutions to PDEs associated with a Sturm-Liouville operator is established in \cite{karafyllis2017siam} by the method of eigenfunction expansion and a finite-difference scheme. We established in this note the ISS in $L^\infty$-norm for a similar setting with considerably simpler computations by De~Giorgi iteration and Lyapunov method. Moreover, the ISS in $L^2$-norm for solutions to certain semilinear parabolic PDEs with Neumann or Robin boundary disturbances is established in \cite{Zheng:2017} by Lyapunov method. These achievements show that the techniques and tools developed in this note and \cite{Zheng:2017} are effective for the application of Lyapunov method to the analysis of the ISS for certain linear and nonlinear PDEs with different type of boundary disturbances.
\end{remark}
\begin{remark}
The method developed in this work can be also applied to linear problems with multidimensional spatial variables, e.g.,
\begin{subequations}
 \begin{align*}
 &u_t-\mu \Delta u+c(x,t)u=f(x,t),\ \ \text{in }\  \Omega\times \mathbb{R}_{+},\\
 &u(x,t)=0\ \ \text{on }\ \Gamma_0,\ u(x,t)=d(t)\ \ \text{on }\ \Gamma_1,\label{}\\
  &u(x,0)=u_0(x),\ \ \text{in }\  \Omega,
\end{align*}
\end{subequations}
where $\Omega\subset\mathbb{R}^n (n\geq 1)$ is an open bounded domain with smooth boundary $\partial \Omega=\Gamma_0\cup\Gamma_1$, $\Gamma_0\cap\Gamma_1=\emptyset$, $ c(x,t)$ is a smooth function in $ \Omega\times \mathbb{R}_{\geq 0}$ with $0< m\leq c(x,t)\leq M$, $\Delta$ is the Laplace operator, and $\mu>0$ is a constant.
Under appropriate assumptions on $\mu,m,M$ and by the technique of splitting and De~Giorgi iteration, it can be shown that the following estimates hold:
\begin{align*}
\|u(\cdot,t)\|_{L^2(\Omega)}\leq & C_0\|u_0\|_{L^2(\Omega)} e^{-\lambda t}
+ C_1\max\limits_{s\in [0,t]}|d(s)| \notag\\
&+ C_2 \max\limits_{(x,s)\in\overline{Q}_t}|f(x,s)|,
\end{align*}
and \begin{align*}
 \|u(\cdot,t)\|_{L^2(\Omega)}^2 \leq& C_0\|u_0\|_{L^2(\Omega)}^2 e^{-\lambda t}+ C_1\max\limits_{s\in [0,t]}|d(s)|^2 \nonumber \\
                     &+ C_2 \int_{0}^t\|f(\cdot,s)\|^2\text{d}s,
\end{align*}
where $\overline{Q}_t=\overline{\Omega}\times [0,t]$, $C_0$, $C_1$, $C_2$, and $\lambda$ are {}{some positive constants independent of $t$}.
 \end{remark}

\section{Concluding Remarks}\label{Sec: Conclusion}
This work applied the technique of De~Giogi iteration to the establishment of ISS properties for nonlinear PDEs. The ISS estimates in $L^2$-norm w.r.t. boundary and in-domain disturbances for Burgers' equation with Dirichlet boundary conditions have been obtained. {The considered setting is a complement of the problems dealt with in \cite{Zheng:2017}, where the ISS in $L^2$-norm has been established for some {semilinear PDEs} with Robin (or Neumann) boundary conditions. It is worth pointing out that the method developed in this note can be generalized for some problems on multidimensional spacial domain and for dealing with ISS properties of PDEs while considering weak solutions (see, e.g., \cite[Ch.~4]{Wu2006}). Finally, as the method of De~Giogi iteration is a well-established tool for regularity analysis of PDEs, we can expect that the method developed in this work is applicable in the study of a wider class of nonlinear PDEs, such as {Chaffee-Infante equation, Fisher-Kolmogorov equation, generalized Burgers' equation, Kuramoto-Sivashinsky equation, and linear or nonlinear Schr\"{o}dinger equations}.}



%
%

\ifCLASSOPTIONcaptionsoff
  \newpage
\fi

\end{document}